\documentclass[preprint,11pt]{elsarticle}
\usepackage{amsfonts}
\usepackage{}
\usepackage{amsmath,amssymb,amsthm}
\usepackage{graphicx}
\usepackage{color}
\usepackage{soul} 
\usepackage{color, xcolor} 

\makeatletter
\def\ps@pprintTitle{%
  \let\@oddhead\@empty
  \let\@evenhead\@empty
  \def\@oddfoot{\reset@font\hfil\thepage\hfil}
  \let\@evenfoot\@oddfoot}
\makeatother

\newcommand{\Q}{\mathbb Q}
\newcommand{\R}{\mathbb{R}}

\newcommand{\N}{\mathbb{N}}
\newcommand{\Z}{\mathbb{Z}}

\newcommand{\lra}{\longrightarrow}
\newcommand{\UU}{\mathcal U}

\newcommand{\dd}{\boldsymbol\delta}

\newtheorem{defin}{Definition}[section]
\newenvironment{definition}{\begin{defin}\rm}{\end{defin}}
\newtheorem{theorem}[defin]{Theorem}
\newtheorem{lemma}[defin]{Lemma}
\newtheorem{proposition}[defin]{Proposition}
\newtheorem{corollary}[defin]{Corollary}
\newtheorem{qu}{Question}
\newtheorem{exa}[defin]{Example}

\begin{document}

\begin{frontmatter}

\title{Uniformization of topological MV-algebras\corref{cor2}}


\author[mymainaddress]{Li-Hong Xie}
\ead{yunli198282@126.com}

\author[mymainaddress1]{Jiang Yang\corref{cor1}} 
\ead{yangjiangdy@126.com}

\cortext[cor1]{Corresponding author: Jiang Yang}
\cortext[cor2]{The project is supported by the Natural Science Foundation of Guangdong Province under Grant (No. 2021A1515010381) and the Innovation Project of Department of Education of Guangdong Province (No. 2022KTSCX145).}
\address[mymainaddress]{School of Mathematics and Computational Science, Wuyi University, Jiangmen, Guangdong, 529000, P.R. China}

\address[mymainaddress1]{ School of Mathematical Sciences, Guangxi Minzu University, Nanning, 530006, P.R. China}

\begin{abstract}
In this paper, First, we prove that for every topological MV-algebra the topology is induced by a uniformity. As a application, we obtain that every topological MV-algebra is completely regular and every \(T_0\) first countable topological MV-algebra is metrizable. Second, we study the completion of the canonical difference uniformity in topological MV-algebra. Finally, we construct some examples show that not all topological MV-algebras are canonically completable and the property of canonically
completable topological MV-algebra does
not imply local order convexity.  
\end{abstract}

\begin{keyword}
uniformity; topological MV-algebra; completely regular; metrizability

\MSC[2020] 06D35; 54E15; 54E35; 54D15
\end{keyword}

\end{frontmatter}

\section{Introduction}
MV-algebras were introduced by Chang \cite{Chang1958} as an algebraic
framework for infinite-valued \L ukasiewicz logic. They extend Boolean
algebras by allowing truth values to vary beyond the classical two-valued
setting. The development of the subject revealed a close connection between
many-valued logic and ordered algebra. In particular, Mundici's
representation theorem \cite{Mundici1986} establishes a categorical
equivalence between MV-algebras and abelian lattice-ordered groups with a
strong unit. Under this equivalence, an MV-algebra is represented as an
interval $\Gamma(G,u)=[0,u]$ in such a group, with
$x\oplus y=(x+y)\wedge u$ and $\neg x=u-x$.
This correspondence provides a structural foundation for the subject;
a systematic account of MV-algebras and their logical and order-theoretic
aspects can be found in \cite{CDM2000}.

A natural question is how this algebraic structure constrains a topology
for which the MV-operations are continuous.
More generally, topological MV-algebras bring together algebraic
operations, lattice order, and convergence, and make it possible to study
the extent to which familiar principles of topological algebra persist
in a setting governed by truncated addition. Hoo initiated the systematic study of topological MV-algebras in
\cite{Hoo1997}. His work investigated topologies constructed from ideals,
including $I$-adic topologies and their completions, and established
basic interactions between algebraic and topological properties. An
important separation result from this early theory is that the axioms
$T_0$, $T_1$, and $T_2$ are equivalent for topological MV-algebras.
Subsequent work developed the relationships among neighborhood systems,
order, and uniform structures. Weber \cite{Weber2012} proved that the
topology of a topological MV-algebra is determined by its neighborhoods
of $0$, and studied the relation between local convexity and uniform
MV-algebras.

The study of filter topologies further clarified the role of algebraic
compatibility. Luan, Weber, and Yang \cite{LWY2021} characterized when a
filter topology turns an MV-algebra into a topological MV-algebra and,
under those conditions, related it to a uniform MV-algebra. Their
results also show why continuity requirements must be stated carefully:
a filter topology need not make the binary MV-operation jointly
continuous. These investigations provide substantial information about
particular constructions and about topologies with additional order
properties. They also lead to a broader question concerning an
arbitrary topology already known to make the MV-operations continuous:
must the underlying topological space admit a compatible uniformity?

There is a strong precedent for a positive answer in the theory of
topological groups. The group operations canonically generate left and
right uniformities, and the resulting uniform structures explain the
complete regularity of topological groups. Moreover, the
Birkhoff--Kakutani metrization theorem asserts that a Hausdorff
topological group is metrizable precisely when it is first countable \cite{AT}. For MV-algebras, however, these conclusions cannot
simply be transferred through the algebraic interval representation.
The operation $\oplus$ is truncated addition, and its translations are
generally not invertible. In addition,
Mundici's algebraic equivalence does not by itself supply a compatible
topological group structure for a given topology on the interval.
Uniformizability and metrization therefore require a topological
argument that takes account of the MV-structure.

Recent results have brought the separation problem into sharper focus.
Gan, Luan, Deng, and Yang \cite{Gan2024} proved that every topological
MV-algebra is regular. Together with the earlier separation theorem,
this yields the equivalence of $T_0$, $T_1$, $T_2$, and $T_3$ in this
class. Regularity, however, does not in general imply complete
regularity, so the existence of enough continuous real-valued functions
to separate points from closed sets requires a further argument.
Li and Yang \cite{LiYang2026} subsequently proved that every
locally convex topological MV-algebra is uniformizable and hence
completely regular. Their theorem establishes the desired conclusion
under an additional order-theoretic hypothesis and makes it natural to
ask whether that hypothesis is needed for the underlying space to be
uniformizable.

The purpose of this paper is to establish uniformizability on every topological MV-algebra and to obtain the corresponding separation and metrization
theorems. The paper organize as follows: In Section \ref{Sec:2}, we establish some basic facts. Section \ref{Sec:3} prove that for every topological MV-algebra the topology is induced by a uniformity. As a application, we obtain that every topological MV-algebra is completely regular and every \(T_0\) first countable topological MV-algebra is metrizable. In Section \ref{Sec:4}, we study the completion of the canonical difference uniformity in topological MV-algebra. We mainly prove that a Hausdorff topological MV-algebra \(A\) is canonically completable if and only if the operation \(\oplus:(A,\mathcal {U}_A)^2\to(A,\mathcal {U}_A)\) is Cauchy continuous, where the domain carries the product uniformity; when this holds, the topological MV-algebra structure on
\(\widehat A\) extending that of \(A\) is unique; its canonical
difference uniformity equals the completion uniformity
\(\widehat\UU_A\), so \(\widehat A\) is canonically complete (see Theorem \ref{thm:criterion}). Section \ref{Sec:5} constructs some examples show that not all topological MV-algebras are canonically completable and the property of canonically
completable topological MV-algebra does
not imply local order convexity.


\section{Preliminaries}\label{Sec:2}

\begin{definition}\label{Def:1}\cite{Chang1958, CDM2000}
An MV-algebra is an algebra \((A,\oplus,^\ast, 0)\), where \((A,\oplus, 0)\) is an abelian monoid and
\(^\ast\) is a unary operation on \( A\) such that for all \(x,y\in A\),
\begin{enumerate}
\item[(MV1)] \(x^{\ast\ast}=x\);
\item[(MV2)] \(x\oplus 0^\ast=0^\ast\);
\item[(MV3)] \((x^\ast\oplus y)^\ast\oplus y=(y^\ast\oplus x)^\ast\oplus x\).
\end{enumerate}
\end{definition}
Given an MV-algebra \((A,\oplus,^*,0)\), for each \(x,y\in A\) we define:
\begin{equation}\label{q1}
 1:=0^*, \quad
 x\odot y:=(x^*\oplus y^*)^*,\quad
 x\ominus y:=x\odot y^*.\quad
\end{equation}
For any \(x,y\in A\), we define \(x \leq y\) if and only if \(x\ominus y = 0\). We call this order a {\bf natural order} on \(A\).

\begin{proposition}\cite[Proposition 1.1.5 ]{CDM2000}\label{pro1}
On each MV-algebra \(A\) the natural order
determines a lattice structure. Specifically, the join \(x\vee y\) and the meet \(x\wedge y\)
of the elements \(x\) and \(y\) are given by
\begin{enumerate}
\item[(a)] \(x\vee y=(x\odot y^*)\oplus y=(x\ominus y)\oplus y\);
\item[(b)] \(x\wedge y=(x^*\vee y^*)^*=x\odot(x^*\oplus y).\)
\end{enumerate}
\end{proposition}

\begin{exa}
Let \([0,1] = \{x \in \mathbb{R} \mid 0 \leq x \leq 1\}\) and for each \(x,y\in [0,1]\) define by:
\[x \oplus y := \min\{1, x + y\}, \qquad
x^{\ast} := 1 - x.\]
Then $([0,1], \oplus, ^\ast, 0)$ is an MV-algebra, known as the {\bf standard MV-algebra}.
\end{exa}

\begin{lemma}\label{Lem:1}
Let $A$ be an \(MV\)-algebra. For any \(x,y,z\in A\), the following statements hold:

\begin{enumerate}
 \item[(a)]  \(x\odot y\le z\) if and only if \(x\le y^*\oplus z\) \cite[Lemma 1.1.4]{CDM2000};
 \item[(b)] \(x\odot (y\vee z)=(x\odot y)\vee(x\odot z)\) \cite[Proposition 1.1.6]{CDM2000};
 \item[(c)] if \(x\leq y\), then \(y^*\leq x^*\) and \(x\odot z\leq y\odot z\) \cite[Lemma 1.1.4]{CDM2000};
 \item [(d)] \((x^*\oplus y)\vee(y^*\oplus x)=1\) \cite[Proposition 1.1.7]{CDM2000};
 \end{enumerate}
\end{lemma}

\begin{lemma}\label{Lem:2}
Let $A$ be an \(MV\)-algebra. For any \(x,y,z\in A\), the following statements hold:
\begin{enumerate}
\item[(a)] \(x\wedge y=x\ominus(x\ominus y);\)
\item[(b)] \((x\ominus y)\ominus z= x\ominus (y\oplus z)=(x\ominus z)\ominus y;\)
 \item[(c)] \((x\wedge y)\ominus z =(x\ominus z)\wedge(y\ominus z)\);
 \item [(d)] \((x\ominus y)\oplus(x\wedge y)=x\);
\item [(e)] \((x\ominus(x\ominus y))\oplus(y\ominus x)=y.\)
 \end{enumerate}
\end{lemma}

\begin{proof}
For (a). By (\ref{q1}) and Proposition \ref{pro1}, we have
\begin{align*}
 x\ominus(x\ominus y)&=x\ominus(x\odot y^*)\\
 &=x\odot(x\odot y^*)^*\\
 &=x\odot(x^*\oplus y)\\
 &=x\wedge y.
\end{align*}

For (b). By (\ref{q1}) and associativity and
commutativity of \(\oplus\), we have
\begin{align*}
 (x\ominus y)\ominus z&=(x\odot y^*)\odot z^*\\
 &=(x^\ast\oplus y)^*\odot z^*\\
 &=((x^\ast\oplus y)\oplus z)^*\\
 &=(x^\ast\oplus (y\oplus z))^*\\
 &=x\odot (y\oplus z)^*\\
 &=x\ominus (y\oplus z)\\
\end{align*}
and
\begin{align*}
 (x\ominus y)\ominus z&=(x\odot y^*)\odot z^*\\
 &=(x^\ast\oplus y)^*\odot z^*\\
 &=((x^\ast\oplus y)\oplus z)^*\\
 &=((x^\ast\oplus z)\oplus y)^*\\
 &=(x\ominus z)\ominus y.
\end{align*}

For (c). Firstly, we shall prove that
\begin{equation}\label{q2}
(b\wedge c) \odot a=(b\odot a)\wedge (c\odot a)
\end{equation}
holds for each \(a,b,c \in A\).

Since \(b\wedge c\leq b\) and \(b\wedge c\leq c\), by Lemma \ref{Lem:1}(c) we have \((b\wedge c)\odot a\leq b \odot a\) and \((b\wedge c)\odot a\leq c \odot a\). Hence,
\((b\wedge c)\odot a\leq (b \odot a)\wedge (c \odot a)\).

Put \(w= (b \odot a)\wedge (c \odot a)\). Then \(w\leq b \odot a\). Hence, by Proposition \ref{pro1}(b), Lemma \ref{Lem:1}(c) and associativity and
commutativity of \(\odot\),
\begin{align*}
w\odot (b^*\oplus c)&\leq (b \odot a)\odot (b^*\oplus c)\\
&=a\odot (b\odot(b^*\oplus c))\\
&=a\odot(b\wedge c)\\
&=(b\wedge c)\odot a
\end{align*}
and likewise \(w\odot (c^*\oplus b)\leq (b\wedge c)\odot a\) by \(w\leq c \odot a\). Therefore, by Proposition \ref{pro1}(b)(d),
\begin{align*}
w&=w\odot((c^*\oplus b)\vee ((b^*\oplus c)))\\
&=(w\odot(c^*\oplus b))\vee(w\odot(b^*\oplus c))\\
&\leq (b\wedge c)\odot a.
\end{align*}
Hence, we have proved that (\ref{q2}): \((b\wedge c) \odot a=(b\odot a)\wedge (c\odot a)\). By (\ref{q1}) and (\ref{q2}), we have that
\begin{align*}
(x\wedge y)\ominus z &=(x\wedge y)\odot z^*\\
&=(x\odot z^*)\wedge(y\odot z^*)\\
&=(x\ominus z)\wedge(y\ominus z).
\end{align*}

For (d). From Proposition \ref{pro1}(a) and Lemma \ref{Lem:2}(a) it follows that
\begin{align*}
(x\ominus y)\oplus(x\wedge y)&=(x\ominus y)\oplus(x\ominus(x\ominus y))\\
&=x\vee (x\ominus y)\\
&=x,
\end{align*}
where the last equality uses \(x\ominus y\le x\).

For (e). This follows from Lemma \ref{Lem:2}(a) and (d).
\end{proof}

Let \(A\) be an MV-algebra. For each \(x,y\in A\), define \(\dd:A\times A\lra A\times A\) by
\begin{equation}\label{eq:directional}
 \dd(x,y)=(x\ominus y,y\ominus x)
\end{equation}
and for \(p=(a,b)\), \(q=(c,d)\) in \(A^2\), define \(\Phi:A^2\times A^2\lra A^2\) by
\begin{equation}\label{eq:F}
 \Phi(p,q)=
 \bigl((a\ominus d)\oplus(c\ominus b),
       (b\ominus c)\oplus(d\ominus a)\bigr).
\end{equation}

\begin{lemma}\label{lem:identities}
Let \(A\) be an MV-algebra. \(\dd\) and \(\Phi \) are defined by (\ref{eq:directional}) and (\ref{eq:F}), respectively. Then
\begin{equation} \label{eq:composition}
 \dd(x,z)=\Phi\bigl(\dd(x,y),\dd(y,z)\bigr)
\end{equation}
hold for each \(x,z, y\in A\).
Moreover, \(\Phi((0,0),(0,0))=(0,0)\).
\end{lemma}

\begin{proof}
Put \(p=x\ominus z, q=y\ominus z\). By Proposition \ref{pro1}(a), Lemma \ref{Lem:2}(b), we have
\begin{align*}
(x\ominus y)\ominus(z\ominus y)&=x\ominus (y\oplus(z\ominus y))\\
&=x\ominus (y\vee z)\\
&=x\ominus (z\oplus(y\ominus z))\\
&=(x\ominus z)\ominus(y\ominus z)=p\ominus q
\end{align*}
and, by Lemma \ref{Lem:2}(a)(c), we have
\begin{align*}
(y\ominus z)\ominus(y\ominus x)&=(y\ominus (y\ominus x))\ominus z\\
&=(y\wedge x)\ominus z\\
&=(y\ominus z)\wedge (x\ominus z)\\
&=p\wedge q.
\end{align*}
Therefore, by Lemma \ref{Lem:2}(d), we have
\[((x\ominus y)\ominus(z\ominus y))\oplus ((y\ominus z)\ominus(y\ominus x))=(p\ominus q)\oplus(p\wedge q)=p=x\ominus z.\]
This proves the first coordinate of \eqref{eq:composition}. Interchanging \(x\) and
\(z\), and using commutativity of \(\oplus\), proves the second coordinate. The assertion about \(\Phi((0,0),(0,0))=(0,0)\) follows immediately from
\eqref{eq:F}.
\end{proof}
\section{Compatible uniformity on every topological MV-algebra}\label{Sec:3}

We begin by introducing some notations. Let $X$ be a set and let $A$ and $B$ be subsets of $X\times X$, i.e. relations on the set $X$. The inverse relation of $A$ will be denoted by $-A$, i.e., $-A=\{(x,y):(y,x)\in A\}$. The composition of $A$ and $B$ will be denoted by $A+B$, i.e., $A+B=\{(x,z): \text{ there is a } y\in X \text{ such that } (x,y)\in A \text{ and } (y,z)\in B\}.$ $2A$ means that $A+A$.

The diagonal of the Cartesian product $X\times X$ is the set $\Delta_X=\{(x,x):x\in X\}$. Every set $V\subseteq X\times X$ that contains $\Delta_X$ and satisfies the condition $V=-V$ is called an {\it entourage of diagonal}. The family of all entourages of diagonal will be denoted by $\mathcal {D}_X$. Let $x_0$ be a point in $X$ and let $V\in \mathcal {D}_X$. The set $V[x_0]=\{x\in X:(x_0,x)\in V\}$ is called the {\it ball with center $x_0$} and {\it radius} $V$.

A {\bf uniformity} \cite{En} on a set $X$ is a subfamily $\mathcal {U}\subseteq \mathcal {D}_X$ which satisfies the following conditions:
\begin{enumerate}
\item[$(U1)$] If $V\in \mathcal {U}$ and $V\subseteq W\in \mathcal {D}_X$, then $W\in \mathcal {U}$.
\item[$(U2)$] If $V_1,V_2\in \mathcal {U}$, then $V_1\cap V_2 \in \mathcal {U}$.
\item[$(U3)$] For every $V\in \mathcal {U} $, there exists a $W\in \mathcal {U}$ such that $2W\subseteq V$.
\end{enumerate}

A uniformity $\mathcal {U}$ on a set $X$ is called a {\bf separated uniformity} if $\mathcal {U}$ satisfies the condition: $(U4)$ $\bigcap \mathcal {U}=\Delta_X.$

A subfamily $\mathcal {B}\subseteq \mathcal {U}$ is called a base for the uniformity $\mathcal {U}$ if for every $V\in \mathcal {U} $ there exists a $W\in\mathcal {B}$ such that $W\subseteq V$.

The following result is well known.
\begin{theorem}\cite[Theorem 8.1.1]{En}\label{t1}
For every uniformity $\mathcal {U}$ on a set $X$ the family $\mathcal {O}=\{U\subseteq X:\text{ for every }x\in U \text{ there is a }V\in\mathcal {U} \text{ such that }V[x]\subseteq U\}$ is a topology on the set.
\end{theorem}

The topology $\mathcal {O}$ is called the {\bf topology induced by the uniformity $\mathcal {U}$}. A topological space $X$ with a topology $\mathcal {T}$ is {\it (pseudo)metrizable} if the topology  $\mathcal {T}$ is induced by a (pseudo)metric $d$ on the set $X$.

 \begin{definition}\cite{Hoo1997}
A topological \(MV\)-algebra is an \(MV\)-algebra \(A\) equipped with a topology such that
\[
        \oplus:A\times A\longrightarrow A,
        \qquad
        {}^\ast:A\longrightarrow A
\]
are continuous. Then the term operations \(\odot,\ominus,\vee,\wedge\) are also continuous.
\end{definition}

Let \(\mathcal N_0\) denote the collection of open neighborhoods of
zero in a topological MV-algebra \(A\).

\begin{theorem}\label{Th}
Let \(A\) be a topological MV-algebra. For \(U\in\mathcal N_0\), put
\begin{equation}\label{eq:EU}
 E_U=\{(x,y)\in A^2:\dd(x,y)\in U\times U\}, \text{~where~} \dd \text{~is defined by~} (\ref{eq:directional}).
\end{equation}
Then \(\{E_U:U\in\mathcal N_0\}\) is a base for the uniformity
\(\UU_A=\{D\in\mathcal {D}_A: E_U\subseteq D \text{~for some~} U\in \mathcal N_0 \text{~such that~} E_U\subseteq D \}\) inducing the given topology of \(A\). If \(A\) is \(T_0\), then \(\UU_A\) is a separated uniformity.
\end{theorem}

\begin{proof}
Firstly, we show that \(\UU_A\) is a uniformity on \(A\).

For \((U1)\). It is trivial.

For \((U2)\). This follows from the fact \(E_U\cap E_V=E_{U\cap V}\) holds for each \(U,V\in \mathcal N_0\).

For \((U3)\). Let \(\Phi:A^2\times A^2\to A^2\) be defined by (\ref{eq:F}), that is,  \[\Phi((a,b),(c,d))=
 \bigl((a\ominus d)\oplus(c\ominus b),
       (b\ominus c)\oplus(d\ominus a)\bigr)\] for each \(((a,b),(c,d))\in A^2\times A^2\). Since \(A\) is a topological MV-algebra, \(\Phi\) is continuous on  \(A^2\times A^2\), which is endowed the the product topology.
       Take any \(D\in \UU_A\). Then there is \(U\in \mathcal N_0\) such that \( E_U\subseteq D\). By Lemma \ref{lem:identities}, there is an open neighborhood \(V\) of \(0\) such that
\[
 \Phi\bigl((V\times V)\times(V\times V)\bigr)
 \subseteq U\times U,
\]
because \(\Phi((0,0),(0,0))=(0,0)\). If \((x,y),(y,z)\in E_V\), identity \eqref{eq:composition} now gives
\(\dd(x,z)\in U\times U\). Hence \(2 E_V\subseteq E_U\),
as required. Clearly, \(\{E_U:U\in\mathcal N_0\}\) is a base for the uniformity
\(\UU_A\).

Secondly, we shall prove that \(\UU_A\) induce the given topology of \(A\). In fact, since \(A\) is a topological MV-algebra, for each fix \(a\in A\), the left translation \(\rho_a:A\rightarrow A\) defined by \(\rho_a(x)=a\ominus x\) for each \(x\in A\) is continuous. Similarly, the right translation \(\rho_a:A\rightarrow A\) defined by \(\lambda_a(x)=x\ominus a\) for each \(x\in A\) is continuous. Thus,
\[
 E_U[a]=\{b\in A:a\ominus b\in U,\ b\ominus a\in U\}=\rho_a^{-1}(U)\cap \lambda_a^{-1}(U).
\]
is an open neighborhood of \(a\) in the original topology for each \(U\in \mathcal N_0\). This implies that the original topology on \(A\) is finer than the topology induced by \(\UU_A\).

Conversely, let \(O\) be an original open neighborhood of \(a\). Define \[\vartheta_a:A^2\rightarrow A\] by \(\vartheta_a(x,y)=(a\ominus x)\oplus y\) for each \((x,y)\in A^2\).
Clearly, \(\vartheta_a\) is continuous and \(\vartheta_a(0,0)=a\). Choose
\(U\in\mathcal N_0\) with \(\vartheta_a(U\times U)\subseteq O\).
For \(x\in E_U[a]\), Lemma \ref{Lem:2}(e) yields
\[
 x=(a\ominus(a\ominus x))\oplus(x\ominus a)=\vartheta_a\bigl(a\ominus x, x\ominus a\bigr)\in O.
\]
Thus \(E_U[a]\subseteq O\). The sets \(E_U[a]\) are consequently a
neighborhood base for the original topology at every \(a\). This implies that the topology induced by \(\UU_A\) is finer than the original topology on \(A\).
Thus we have proved that \(\UU_A\) induces the given topology of \(A\).

If \(A\) is \(T_0\), then \(A\) is Hausdorff by \cite[Corollary 3.10]{Gan2024}. Hence,
\(\bigcap_{U\in\mathcal N_0}U=\{0\}\), so membership in every \(E_U\)
forces both directional differences to be zero and hence \(x=y\). This implies that the condition: $(U4)$ $\bigcap \UU_A=\Delta_A$ holds.
\end{proof}

\begin{definition}
A topological space \(X\) is called
\begin{enumerate}
\item [(a)]  a  \(T_0\) space if for every pair of distinct points \(x,y\in X\), there exists an open set \(U\in\tau\) such that exactly one of
\(x\) and \(y\) belongs to \(U\);
\item [(b)] a \(T_1\)-space if for any distinct points \(x, y \in X\) the point \(x\) has a neighborhood \(U_x \subseteq X\) such that \(y \notin U_x\);
\item [(c)] a \(T_2\)-space if for any distinct points \(x, y \in X\) the point \(x\) has a neighborhood \(U_x \subseteq X\) such that \(y \notin \overline{U_x}\);
\item [(d)] regular if for any point \(x \in X\) and a neighborhood \(O_x \subseteq X\) of \(x\) there is a neighborhood
\(U_x \subseteq X\) of \(x\) such that \(\overline{U_x}\subseteq O_x\);
\item [(e)] a \(T_3\)-space if \(X\) is a regular \(T_1\)-space;
\item [(f)] completely regular if for every closed set \(F\subseteq X\) and every point \(x\in X\setminus F\), there exists a continuous function
\(f:X\to[0,1]\)
such that \(f(x)=0 \) and \(f(y)=1\) for all \( y\in F;\)
\item [(g)] a Tychonoff space is a completely regular \(T_1\) space.
\end{enumerate}
\end{definition}

\begin{corollary}
Every topological MV-algebra is completely regular.
\end{corollary}

\begin{proof}
Since every topological space the topology of which is induced by a uniformity is completely regular \cite{En}, the statement follows from Theorem \ref{Th}.
\end{proof}
\begin{corollary}
Let \(A\) be a topological MV-algebra. Then the following statements are equivalent:
\begin{enumerate}
\item [(a)] \(A\) is a \(T_0\) space;
\item [(b)] \(A\) is a \(T_1\) space;
\item [(c)] \(A\) is a \(T_2\) space;
\item [(d)] \(A\) is a \(T_3\) space;
\item [(e)] \(A\) is a Tychonoff space.
\end{enumerate}
\end{corollary}

\begin{proof}
The equivalences of \((a)\)-\((d)\) follows from \cite[Corollary 3.10]{Gan2024}. Since every \(T_3\) topological space the topology of which is induced by a uniformity is Tychonoff space\cite[Corollary 8.1.13]{En}, the statement follows from Theorem \ref{Th}.
\end{proof}

\begin{lemma}\label{Lem:3}\cite[Theorem 8.1.21]{En}
A uniformity \(\UU\) on a set \(X\) is induced by a metric on \(X\) if and only if \(\UU\) is separated with a countable base.
\end{lemma}

\begin{corollary}\label{Co1}
Every \(T_0\), first countable topological MV-algebra is metrizable.
\end{corollary}

\begin{proof}
This follows from Lemma \ref{Lem:3} and Theorem \ref{Th}.
\end{proof}

\section{The completion of the canonical difference uniformity in topological MV-algebra}\label{Sec:4}
In this section, we assume that every topological MV-algebra is Hausdorff. In this section, we study the Hausdorff completion of the compatible uniformity on a
topological MV-algebra defined in Theorem \ref{Th}. We call this uniformity the {\bf canonical difference uniformity} on a topological MV-algebra. Denote the canonical
difference uniformity on a topological MV-algebra \(A\) by \(\mathcal {U}_A\).

Let \((X,\mathcal {U})\) be a uniform space. A family \(\mathcal {F}\) of subsets in \(X\) is {\bf Cauchy}
if, for every entourage \(E\in \mathcal {U}\), some \(F\in\mathcal {F}\) satisfies \(F\times F\subseteq E\). A filter (resp. a filter basis)\(\mathcal {F}\) is called a {\bf Cauchy filter} (resp. a {\bf Cauchy filter basis} )if \(\mathcal {F}\) is Cauchy \cite[0.34 Definition]{RD1981}. A uniform space \((X,\mathcal {U})\) is called {\bf complete} if every \(\mathcal {U}\)-Cauchy filter on
\((X,\mathcal {U})\) converges to some \(x\) in \((X,\mathcal {U})\). Let both \((X,\mathcal {U}_X)\) and (\(Y,\mathcal {U}_Y\)) be uniform spaces. A map \(f:(X,\mathcal {U}_X)\rightarrow (Y,\mathcal {U}_Y)\) is called {\bf uniformly continuous} if for every \(M\in \mathcal {U}_Y\) there exists \(N\in \mathcal {U}_X\) such that \((f(x),f(y))\in M\) for each \((x,y)\in N\), where \(f\) is called a {\bf uniform equivalence} if \(f\) is bijective and both \(f:(X,\mathcal {U}_X)\rightarrow (Y,\mathcal {U}_Y)\) and its inverse \(f^{-1}:(Y,\mathcal {U}_Y)\rightarrow (X,\mathcal {U}_X)\) are uniformly continuous \cite[0.14 Definition]{RD1981}. A map \(f:(X,\mathcal {U}_X)\rightarrow (Y,\mathcal {U}_Y)\) is called a {\bf uniform embedding} if the induced map \((X,\mathcal {U}_X)\rightarrow (f(X), \mathcal {U}_Y|f(X)), x\mapsto f(x)\) is a uniform equivalence \cite[0.43 Definition]{RD1981}. A map \(f:(X,\mathcal {U}_X)\rightarrow (Y,\mathcal {U}_Y)\) is {\bf Cauchy continuous} if it sends Cauchy filters
to Cauchy filter basises.  A {\bf completion} of a uniform space \((X,\mathcal {U}_X)\) is a pair \(((Y,\mathcal {U}_Y),i)\) consisting of a complete uniform space \((Y,\mathcal {U}_Y)\) and a uniform embedding \(i: (X,\mathcal {U}_X)\rightarrow (Y,\mathcal {U}_Y)\) such that \(i(X)\) is dense in \(Y\) \cite[0.44 Definition]{RD1981}. For each Hausdorff uniform space \((X,\mathcal {U}_X)\) there is a unique Hausdorff completion \(((\widehat{X},\widehat{\mathcal {U}_X}),i)\) of \((X,\mathcal {U}_X)\) in the sense of uniform equivalence \cite[0.48 Theorem, 0.51 Proposition]{RD1981}. In this section, the term ``completion of \((X,\mathcal {U}_X)\) " will then also be applied to the space \((\widehat{X},\widehat{\mathcal {U}_X})\) and, still more abusively, to the set \(\widehat{X}\).

\begin{lemma}\label{lem:extension}
Let \(X\) be dense in a Hausdorff uniform space \(Z\), let \(Y\) be a complete Hausdorff uniform space,
and let \(f:X\to Y\) be continuous. Then \(f\) has a unique continuous extension \(\widehat f:Z\to Y\) if \(f\) is Cauchy continuous.
\end{lemma}

\begin{proof}
Assume \(f\) is Cauchy continuous. For \(p\in Z\), the traces on \(X\) of
neighborhoods of \(p\) generate a Cauchy filter \(\mathcal {F}_p\) in \(X\).
The image \(f(\mathcal {F}_p)\) is a Cauchy filter basis, and therefore has a unique limit
in \(Y\), because \(Y\) is a complete Hausdorff uniform space. Define
\[
 \widehat f(p)=\lim f(\mathcal {F}_p).
\]
For \(p\in X\), continuity of \(f\) shows that this is \(f(p)\).

To check the continuity of \(\widehat f\), let \(W\) be an open neighborhood of
\(\widehat f(p)\). Choose a closed neighborhood \(C\) of
\(\widehat f(p)\) contained in \(W\), using the regularity of \(Y\).
There is an open neighborhood \(O\) of \(p\) in \(Z\)
with \(f(O\cap X)\subseteq C\).
For every \(q\in O\), the filter basis \(f(\mathcal {F}_q)\) is eventually in
\(C\), so its limit belongs to \(C\). Hence
\(\widehat f(O)\subseteq C\subseteq W\).
Uniqueness follows because \(X\) is dense and \(Y\) is Hausdorff.
\end{proof}

\begin{proposition}\label{prop:homomorphisms}
A continuous MV-homomorphism \(f:A\to B\) is uniformly continuous
from \((A,\UU_A)\) to \((B,\UU_B)\). If \(A\) is an MV-subalgebra
of \(B\) with the subspace topology, then
\(\UU_A=\UU_B|_{A\times A}\).
\end{proposition}

\begin{proof}
Given a zero-neighborhood \(U\) in \(B\), choose a zero-neighborhood
\(V\) in \(A\) with \(f(V)\subseteq U\). Since \(f\) preserves
directional differences, \((f\times f)(E_V)\subseteq E_U\).
For a subalgebra, the identities
\(E_U^B\cap A^2=E_{U\cap A}^A\) and the description of subspace
zero-neighborhoods prove the second assertion.
\end{proof}

\begin{lemma}
\label{lem:dense-uniformity}
Suppose two uniformities on a set \(Y\) induce the same topology
and have the same restriction to a dense subset \(X\).
Then the two uniformities are equal.
\end{lemma}

\begin{proof}
Let \(E\) be a closed entourage for the first uniformity.
Equality of the restrictions gives an entourage \(F\) of the second
uniformity with \(F\cap X^2\subseteq E\cap X^2\).
Choose a product-topologically open entourage \(G\subseteq F\)
for the second uniformity. Such entourages form a base, and the
product topology is the same for both uniformities.
Density of \(X^2\) gives
\[
 G\subseteq\overline{G\cap X^2}\subseteq E.
\]
Thus \(E\) is also an entourage of the second uniformity.
Closed entourages form a base for every uniformity, so the first
uniformity is contained in the second. Interchanging them proves
the reverse containment.
\end{proof}

\begin{definition}\label{De1}
A Hausdorff topological MV-algebra \(A\) is {\bf canonically complete}
if its canonical difference uniformity is complete.
It is {\bf canonically completable} if the completion \(\widehat{A}\) of the uniform space \((A,\mathcal {U}_A)\) is a topological MV-algebra containing \(A\) as an MV-subalgebra.
\end{definition}

Here \(\widehat A\) carries its uniform completion topology, and
the containment is understood through the canonical uniform embedding
\(j_A:A\to\widehat A\): the operations on \(\widehat A\) must
restrict to those of \(A\).

From now on a topological MV-algebrav \(A\) is Hausdorff and \(\widehat A\) means the
completion of \((A,\UU_A)\), with completion uniformity
\(\widehat\UU_A\).
\begin{theorem}\label{thm:criterion}
A Hausdorff topological MV-algebra \(A\) is canonically completable if and only if the operation \(\oplus:(A,\mathcal {U}_A)^2\to(A,\mathcal {U}_A)\) is Cauchy continuous, where the domain carries the product uniformity.
When this holds, the topological MV-algebra structure on
\(\widehat A\) extending that of \(A\) is unique. Its canonical
difference uniformity equals the completion uniformity
\(\widehat\UU_A\), so \(\widehat A\) is canonically complete.
\end{theorem}

\begin{proof}
Identify \(A\) with its image under \(j_A\), and write
\(\widehat\UU_A\) for the uniformity of its Hausdorff completion.
Thus \(A\) is a dense uniform subspace of \(\widehat A\).

\emph{Necessity.}
Suppose that \(A\) is canonically completable, and denote the
continuous addition on its completion by \(\widehat\oplus\).
Let \(\mathcal {F}\) be a Cauchy filter on \((A,\UU_A)^2\).
The product \((\widehat A,\widehat\UU_A)^2\) is complete. The image of \(\mathcal {F}\) under the uniform
embedding \[j_A \times j_A: A\times A\hookrightarrow\widehat A\times \widehat A\] therefore converges
to some \((p,q)\in\widehat A\times \widehat A\).
Joint continuity gives
\[
 \widehat\oplus(\mathcal {F})\longrightarrow p\,\widehat\oplus\,q
 \qquad\text{in }\widehat A.
\]
Since every convergent family on a uniform space is Cauchy \cite[Proposition 0.35]{RD1981}, \(\widehat\oplus(\mathcal {F})\) is a Cauchy filter basis in \(\widehat A\). Because
\(\widehat\oplus\) restricts to \(\oplus\) and
\(\widehat\UU_A|_{A\times A}=\UU_A\), it follows that
\(\widehat\oplus(\mathcal {F})=\oplus(\mathcal {F})\) is a Cauchy filter basis in \((A,\UU_A)\).
This proves Cauchy continuity.

\emph{Sufficiency.}
Assume that \(\oplus\) is Cauchy continuous.
The identities
\[
 x^*\ominus y^*=y\ominus x,\qquad
 y^*\ominus x^*=x\ominus y
\]
give, for every zero-neighborhood \(U\),
\[
 (x,y)\in E_U\quad\Longleftrightarrow\quad(x^*,y^*)\in E_U.
\]
Thus the involution is a uniform isomorphism of \((A,\UU_A)\)
onto itself. The extension theorem for uniformly continuous maps
into complete Hausdorff uniform spaces gives a unique uniformly
continuous map \({}^{\widehat{*}}:\widehat A\to\widehat A\) with
\({}^{\widehat{*}}|_A={}^*\). Both \({}^{\widehat{*}}\circ{}^{\widehat{*}}\) and the identity map are
continuous and agree on the dense subset \(A\); hence
\({}^{\widehat{*}}\circ{}^{\widehat{*}}=\mathrm{id}_{\widehat A}\).

Assume  \(\oplus:(A,\mathcal {U}_A)^2\to(A, \mathcal {U}_A)\) is Cauchy continuous. Then so is \(\oplus:(A,\mathcal {U}_A)^2\to(\widehat{A},\widehat{\mathcal {U}}_A)\), because \(\widehat\UU_A|_{A\times A}=\UU_A\). Since \(A\times A\) is dense in the uniform space \(\widehat{A}\times \widehat{A}\) and \(\widehat{A}\) is a complete Hausdorff uniform space, applying Lemma \ref{lem:extension} to \(\oplus: A\times A\rightarrow \widehat{A}\) there is a unique continuous extension \(\widehat{\oplus}: \widehat{A}\times \widehat{A}\rightarrow \widehat{A}\) of \(\oplus\).

Put \(\widehat0=j_A(0)\). Every term formed from
\(\widehat\oplus,{}^{\widehat{*}},\widehat0\) defines a continuous map on
the appropriate finite power of \(\widehat A\).
For an MV-identity in \(n\) variables, its two term maps agree
on \(A^n\), because the operations restrict to those of \(A\).
The set on which two continuous maps into a Hausdorff space agree
is closed. Since \(A^n\) is dense in \(\widehat A^n\),
the identity holds throughout \(\widehat A^n\).
In particular, for \(p,q,r\in\widehat A\), we obtain
\begin{align*}
 p\,\widehat\oplus\,q&=q\,\widehat\oplus\,p,\\
 (p\,\widehat\oplus\,q)\,\widehat\oplus\,r
   &=p\,\widehat\oplus\,(q\,\widehat\oplus\,r),\\
 p\,\widehat\oplus\,\widehat0&=p,
       \qquad p^{\widehat{*}\widehat{*}}=p,\\
 p\,\widehat\oplus\,\widehat0^{\widehat{*}}&=\widehat0^{\widehat{*}},\\
 (p^{\widehat{*}}\,\widehat\oplus\,q)^{\widehat{*}}\,\widehat\oplus\,q
   &=(q^{\widehat{*}}\,\widehat\oplus\,p)^{\widehat{*}}\,\widehat\oplus\,p.
\end{align*}
These are all the defining MV-axioms. Consequently
\((\widehat A,\widehat\oplus,{}^{\widehat{*}},\widehat0)\) is a Hausdorff
topological MV-algebra, and \(j_A\) is an MV-embedding.
Definition~\ref{De1} now shows that
\(A\) is canonically completable.

Two continuous extensions of addition agree on the dense set
\(A^2\), and two continuous extensions of the involution agree
on the dense set \(A\). Since \(\widehat A\) is Hausdorff,
the respective extensions agree everywhere; the zero is forced
to be \(j_A(0)\). This proves uniqueness.

Let \(\mathcal V\) be the canonical difference uniformity of
the resulting topological MV-algebra \(\widehat A\).
By Theorem~\ref{Th}, it induces the completion
topology. Since \(A\) is a topological MV-subalgebra,
Proposition~\ref{prop:homomorphisms} gives
\[
 \mathcal V|_{A^2}=\UU_A=\widehat\UU_A|_{A^2}.
\]
The two uniformities therefore induce the same topology and have
the same restriction to the dense subset \(A\).
Lemma~\ref{lem:dense-uniformity} yields
\(\mathcal V=\widehat\UU_A\). As the latter is complete,
\(\widehat A\) is canonically complete.
\end{proof}

\begin{corollary}
\label{cor:universal}
Let \(A\) be a canonically completable topological MV-algebra. Every continuous
MV-homomorphism \(f:A\longrightarrow B\) into a canonically complete Hausdorff
topological MV-algebra extends uniquely to a continuous
MV-homomorphism \(\widehat f:\widehat A\longrightarrow B\).
Moreover, \(A\) is canonically completable if and only if it
admits a topological MV-embedding into some canonically complete
Hausdorff topological MV-algebra.
\end{corollary}

\begin{proof}
By Theorem~\ref{thm:criterion}, \(\widehat A\)
is a canonically complete topological MV-algebra extending \(A\).
By Proposition~\ref{prop:homomorphisms}, \(f\) is uniformly
continuous, so it extends uniquely to a uniformly continuous
map \(\widehat f:\widehat A\longrightarrow B\) by \cite[Theorem 8.3.10]{En}.
Preservation of the operations follows by continuity and density.

One direction of the embedding assertion follows by taking
\(\widehat A\). Conversely, identify \(A\) with a topological
MV-subalgebra of a canonically complete algebra \(B\).
Its closure \(C\) in \(B\) is a subalgebra \cite[Theorem 3.5]{Hoo1997}. It is a closed uniform subspace
of a complete space \cite[Proposition 0.39]{RD1981}, hence complete.
By Proposition~\ref{prop:homomorphisms}, its canonical uniformity
is the subspace uniformity from \(B\), and \(A\) is uniformly
dense in \(C\). Thus \(C\) realizes the uniform completion of
\(A\) and supplies the required continuous operations.
\end{proof}

According to \cite{Na1965}, a subset \(V\) of a partially ordered set is {\bf order convex} if
\(a,b\in V\) and \(a\le x\le b\) imply \(x\in V\). A topology is
{\bf locally order convex} if every point has a neighborhood base
consisting of order convex sets. Let \(A\) be a MV-algebra and \(V\) is a subset of \(A\). \(V\) is called a {\bf downward closed set} if \(z\leq v\) for some \(v\in V\), then \(z\in V\).
\begin{theorem}\label{thm:uniform}
For a Hausdorff topological MV-algebra \(A\), the following are
equivalent:
\begin{enumerate}
\item[(a)] \(A\) is locally order convex;
\item[(b)] zero has a neighborhood base of downward closed sets;
\item[(c)] \(\oplus:(A,\UU_A)^2\longrightarrow(A,\UU_A)\) is uniformly continuous, where the domain carries the product uniformity;
\item[(d)] the MV-operations extend to uniformly continuous operations
  on \((\widehat A,\widehat\UU_A)\).
\end{enumerate}
When this holds, \(A\) is canonically completable and its completion \(\widehat{A}\) is again locally order convex.
\end{theorem}

\begin{proof}
An order convex set containing zero is downward closed, so
(a) implies (b).
Assume (b). For each \(a\in A\) given a neighborhood \(O\) of \(a\), choose an open
zero-neighborhood \(U\) with \(E_U[a]\subseteq O\), using
Theorem~\ref{Th}. Choose a downward closed
zero-neighborhood \(V\subseteq U\). Then \(E_V[a]\) is a
neighborhood of \(a\) and \(E_V[a]\subseteq E_U[a]\subseteq O \).
If \(b,c\in E_V[a]\) and \(b\le z\le c\), using commutativity of \(\odot\), the definition of \(\ominus\) by (\ref{q1}) and Lemma \ref{Lem:1}(c) gives
\[
 a\ominus z\le a\ominus b\in V,\qquad
 z\ominus a\le c\ominus a\in V.
\]
Thus \(z\in E_V[a]\) and \(E_V[a]\) is order convex neighborhood, proving (a).

 (b) \(\Rightarrow\) (c). Take any \(a,b,c,d\in A\). First note the MV-inequality
\begin{equation}\label{eq:sum-estimate}
 (a\oplus c)\ominus(b\oplus d)
 \le(a\ominus b)\oplus(c\ominus d).
\end{equation}
Indeed, \(a\le b\oplus(a\ominus b)\) and
\(c\le d\oplus(c\ominus d)\). Adding and applying
Lemma \ref{Lem:1}(a) proves the inequality.

Given a zero-neighborhood \(U\), choose a downward closed
zero-neighborhood \(W\subseteq U\). Continuity at \((0,0)\)
provides an open zero-neighborhood \(V\) with
\(V\oplus V\subseteq W\).
If \((a,b),(c,d)\in E_V\), then \(a\ominus b,b\ominus a \in V\) and \(c\ominus d,d\ominus c \in V\). Hence, \[(a\ominus b)\oplus (c\ominus d)\in V\oplus V\subseteq W\]
Thus, by \eqref{eq:sum-estimate} and \(W\) being downward closed,
\[(a\oplus c)\ominus(b\oplus d)\in W\subseteq U\]
Similarly, we have \((b\oplus d)\ominus (a\oplus c)\in U\). This implies that \(((a\oplus c),(b\oplus d))\in E_U\). This is uniform continuity.

(c) \(\Rightarrow\) (b). The identities
\[
 x^*\ominus y^*=y\ominus x,\qquad
 y^*\ominus x^*=x\ominus y
\]
give, for every zero-neighborhood \(U\),\[ (x,y)\in E_U\quad\Longleftrightarrow\quad(x^*,y^*)\in E_U.\]
Thus, the involution is always uniformly continuous, so every MV-term, in particular
\(\wedge\), is uniformly continuous.
Given a zero-neighborhood \(U\), choose an open zero-neighborhood
\(V\) such that
\[
 (x,x'),(z,z')\in E_V
 \quad\Longrightarrow\quad
 (x\wedge z,x'\wedge z')\in E_U.
\]
For \(x\in V\), the pair \((x,0)\) lies in \(E_V\); and
\((z,z)\in E_V\) for every \(z\). Therefore \(x\wedge z\in U\), which implies that \(\downarrow V\subseteq U\), where \(\downarrow V=\{y\in A: y\leq x \text{~for some~} x\in V\}\).
Clearly, \(
 V\subseteq\downarrow V
\). Consequently, the downward closed set \(\downarrow V\) is a zero-neighborhood.
This proves (b).

Condition (c), together with uniform continuity of the
involution, gives uniformly continuous extensions of the basic
operations to the completions  by \cite[Theorem 8.3.10]{En}. The MV-identities pass to the
completion as in Theorem~\ref{thm:criterion}; hence (d).
Conversely, restricting a uniformly continuous extension gives
(c). Finally, Theorem~\ref{thm:criterion} identifies the
completion uniformity with the completed algebra's canonical
uniformity. Applying (c) implies (a) to that algebra proves
the last assertion.
\end{proof}

The relationship in Theorem~\ref{thm:uniform} is consistent with
the established theory of lattice and MV-uniformities
\cite{LiYang2026, Weber1991,Weber2012}. It is stronger than the
continuous-extension criterion of Theorem~\ref{thm:criterion}.
Section~\ref{Sec:5} demonstrates that the two criteria
are genuinely different even for a countable Boolean algebra.

Let \((X,\mathcal {U})\) be a uniform space. \((X,\mathcal {U})\) is called {\bf totally bounded} if for each \(U\in \mathcal {U}\), there is a finite set \(F\subseteq X\) such that \(X=\bigcup_{x\in F}U[x]\).
\begin{corollary}\label{cor:precompact}
Let \(A\) be a Hausdorff topological MV-algebra such that the uniform space \((A,\UU_A)\) is totally bounded. Then \(A\) is canonically
completable if and only if \(A\) is locally order convex.
In this case its canonical completion is a compact Hausdorff
topological MV-algebra.
\end{corollary}

\begin{proof}
Assume that \(A\) is locally order convex. Then it follows that \(A\) is canonically
completable from Theorem \ref{thm:uniform}.

Conversely, Since the completion of a totally bounded Hausdorff uniform space is
compact \cite[Corollary 8.3.17]{En}, \(\widehat{A}\) is compact topological MV-algebra. Hence, the continuous extension \(\widehat\oplus:(\widehat{A},\widehat{\mathcal {U}}_A)\times(\widehat{A},\widehat{\mathcal {U}}_A)\longrightarrow (\widehat{A},\widehat{\mathcal {U}}_A)\) of \(\oplus\) is uniformly continuous by \cite[Corollary 0.54]{RD1981}, because its domain
\(\widehat A\times \widehat A\) is compact.
Its restriction \(\oplus\) to \((A,\UU_A)\times (A,\UU_A)\) is therefore uniformly continuous,
and Theorem~\ref{thm:uniform} gives local order convexity.
\end{proof}
\section{Some examples}\label{Sec:5}
In this section, we give two examples which show that not all topological MV-algebras are canonically
completable and the property of canonically
completable topological MV-algebra does
not imply local order convexity.

The following construction shows that first countability, second
countability, and linearity of the MV-order do not suffice for
continuous extendability.

\begin{exa}\label{thm:chain}
There exists a countable Hausdorff second-countable topological
MV-chain \(A\) with canonical completion
\[
 \widehat A=[0,1]\times\Z_{\mathrm{discrete}},
\]
such that even one fixed truncated translation
\(x\mapsto a\oplus x\) has no continuous extension to
\(\widehat A\). In particular, \(A\) is not canonically completable.
\end{exa}

\begin{proof}
Fix an irrational real number \(\alpha\), and put
\[
 G=\Q+\Z\alpha,\qquad A=G\cap[0,1].
\]
Every \(x\in G\) has a unique expression \(x=q+n\alpha\),
where \(q\in\Q\) and \(n\in\Z\).
Define the additive map \(\kappa:G\to\Z\) by
\(\kappa(q+n\alpha)=n\).
The interval \(A\) is an MV-subalgebra of the standard
algebra, with
\[
 x\oplus y=\min\{x+y,1\},\qquad x^*=1-x.
\]
Its natural order is the usual order, so it is an MV-chain.

Define a metric on \(A\) by, for each \(x,y\in A\)
\begin{equation}\label{eq:chain-metric}
d(x,y)=
\begin{cases}
|x-y|, & \kappa(x)=\kappa(y) \\
|x-y|+1,  &\kappa(x)\ne\kappa(y).
\end{cases}
\end{equation}
It is a metric because it is the sum of the usual metric and
the pullback of the discrete metric on \(\Z\).
For each \(n\in\Z\), the set
\[
 A_n=(n\alpha+\Q)\cap[0,1]
\]
is an open and closed subspace with its usual real subspace
topology. Since \(A\) is countable and metrizable, it is
Hausdorff and second countable.

We verify continuity of the MV-operations.
Because \(1\in\Q\),
\(\kappa(1-x)=-\kappa(x)\). One can easily show that the involution \({}^*\) is an
isometry for \eqref{eq:chain-metric}, that is, \(d(x,y)=d(x^*,y^*)\). Hence, \({}^*\) is continuous.
Next, we shall show that \(\oplus\) is cintinuous. Fix \(a,b\in A\) and \(\epsilon >0\).

Case 1: \(a+b<1\). Put \(\delta=\text{min}\{\frac{1-(a+b)}{3}, \frac{\epsilon}{2}\}\). We shall prove that \(d(x\oplus y, a\oplus b)<\epsilon\) whenever \(d(x,a)<\delta\) and \(d(y,b)<\delta\).
Since \(\delta<1\), we have \(\kappa(x)=\kappa(a)\) and \(\kappa(y)=\kappa(b)\). Hence, \(x+y\leq a+b+|x-a|+|y-b|\leq a+b+2\delta<1\). Furthermore,
\begin{align*}
\kappa(x\oplus y)&=\kappa(x+ y)\\
&=\kappa(x)+\kappa(y) \\
&=\kappa(a)+\kappa(b) \\
&=\kappa(a+b)\\
&=\kappa(a\oplus b).
\end{align*}
By (\ref{eq:chain-metric}), we have
\begin{align*}
d(x\oplus y, a\oplus b)&=|(x\oplus y)-(a\oplus b)|\\
&=|(x+ y)-(a\oplus b)|\\
&\leq |x-a|+|y-b|\\
&<2\delta\leq \epsilon
\end{align*}

Case 2: \(a+b>1\). Put \(\delta=\text{min}\{\frac{a+b-1}{3},\frac{1}{2}\}\). Then
\begin{align*}
x+y&\geq a+b-|x-a|-|y-b|\\
&\geq a+b-2\delta>1
\end{align*}
whenever \(|x-a|< \delta\) and \(|y-b|< \delta\). Hence, \[d(x\oplus y, a\oplus b)=d(1,1)=0<\epsilon\] whenever \(d(x,a)<\delta\) and \(d(y,b)<\delta\).

Case 3: \(a+b=1\). Then \(\kappa(a\oplus b)=\kappa(1)=0\). Put \(\delta=\text{min}\{\frac{1}{2},\frac{\epsilon}{3}\}\). Let \(d(x,a)< \delta\) and \(d(y,b)< \delta\). Then \(\kappa(x)=\kappa(a)\) and \(\kappa(y)=\kappa(b)\). Hence, if \(x+y\leq 1\),
then
\begin{align*}
\kappa(x\oplus y)&=\kappa(x+y)\\
&=\kappa(x)+\kappa(y)\\
&=\kappa(a)+\kappa(b)\\
&=\kappa(a+b).
\end{align*}
Thus, by (\ref{eq:chain-metric}), we have \[d(x\oplus y, a\oplus b)=|x+y-(a+b)|\leq |x-a|+|y-b|\leq 2\delta<\epsilon.\]
If \(x+y> 1\), then \(x\oplus y=1\). Hence, \(\kappa(x\oplus y)=\kappa(1)=\kappa (a\oplus b)\). Thus, \[d(x\oplus y, a\oplus b)=d(1,1)=0<\epsilon.\]
This proves continuity of \(\oplus\). Therefore \(A\) with the topology induced by \(d\) defined in (\ref{eq:chain-metric}) is a topological MV-algebra.

A zero-neighborhood base is \(\{U_\epsilon: 0<\varepsilon<1\}\), where
\begin{align*}
 U_\varepsilon&=\{x\in A:d(0,x)<\epsilon\}\\
 &=\{x\in A_0:d(0,x)<\epsilon\}\\
 &=\Q\cap[0,\varepsilon).
\end{align*}
 Put \[ E_{U_\varepsilon}=\{(x,y)\in A\times A:x\ominus y\in U_\epsilon \text{~and~}y\ominus x\in U_\epsilon\}\]
and \[E^\epsilon_d=\{(x,y)\in A\times A:d(x,y)<\epsilon\}.\]

If \((x,y)\in E^\epsilon_d\), then \(d(x,y)<\epsilon<1\). Thus, we have \(\kappa(x)=\kappa(y)\). This implies that \(x-y\in \Q\) and \(|x-y|=d(x,y)<\epsilon\). Since the standard-algebra difference is
\(x\ominus y=\max\{x-y,0\}\), we have \(x\ominus y \in U_\epsilon\) and \(y\ominus x \in U_\epsilon\), proving \(E^\epsilon_d\subseteq  E_{U_\varepsilon}\).

If \((x,y)\in E_{U_\varepsilon}\), then \(x\ominus y \in U_\epsilon\) and \(y\ominus x \in U_\epsilon\). Hence \(\kappa(x\ominus y)=\kappa(y\ominus x)=0\) and \(|x-y|<\epsilon\). Therefore we have \(\kappa(x)=\kappa(y)\). Then \(d(x,y)=|x-y|<\epsilon\), which implies that \((x,y)\in E^\epsilon_d\), proving \(E^\epsilon_d\supseteq  E_{U_\varepsilon}\).
Thus  the canonical difference uniformity \(\UU_A\) on \(A\) is induced by the metric
\eqref{eq:chain-metric}.

Give \([0,1]\times\Z\) the complete metric

\begin{equation}\label{eq:chain-metric2}
\rho((s,n),(t,m))=
\begin{cases}
|s-t|, & n=m \\
|s-t|+1,  &n\ne m.
\end{cases}
\end{equation}
Completeness follows because a Cauchy sequence eventually
is contained in \([0,1]\times\{n\}\) for some \(n\in \Z\) and \([0,1]\times\{n\}\) is complete.
The map \(j:A\rightarrow [0,1]\times\Z\) defined by
\begin{equation}\label{eq:chain-embedding}
 j(x)=(x,\kappa(x))
\end{equation}
is isometric. Since each \(n\alpha+\Q\) is dense in \(\R\), \(A_n=(n\alpha+\Q)\cap [0,1]\) is dense in [0,1] for each \(n\in \Z\).
Hence \(j(A_n)\) is dense in the branch \([0,1]\times\{n\}\) for each \(n\in \Z\). This implies that \(j(A)\) is dense in \([0,1]\times\Z\).
This proves that the uniform space \([0,1]\times\Z\) with the uniformity induced by the metric \(\rho\) defined by (\ref{eq:chain-metric2}) is the canonical completion of the uniform space \((A,\mathcal {U}_A)\).

Now take \(a=1/2\in A_0\). We shall show that the left translation \(\lambda_a:A\rightarrow A\) defined by \(\lambda_a(x)=a\oplus x\) for each \(x\in A\) has no continuous extension to
\(\widehat A\).

 Since \(A_1=(\alpha+\Q)\cap [0,1]\) is dense in [0,1] with the Euclidean topology, we can choose a sequence \(\{b_n:n\in \N\}\subseteq A_1\)
such that \(\{b_n:n\in \N\}\) converges to \(1/2\) in [0,1] and
\[
 b_{2n}<1/2,\quad b_{2n+1}>1/2
\]
for each \(n\in \N\). Hence, the sequence \(\{b_n:n\in \N\}\) is \(d\)-Cauchy in \(A\) and \(
 \{j(b_n):n\in \N\}\) converges to \((1/2,1)\) in \([0,1]\times\Z=\widehat{A}\).

If \(\lambda_a\) has a continuous extension \(\widehat{\lambda_a}\) to \(\widehat A\), then \(\widehat{\lambda_a}(j(b_n))\rightarrow \widehat{\lambda_a}((1/2,1))\). This is impossible.
Since
\[
 a\oplus b_{2n}=1/2+b_{2n},\qquad
 \kappa(a\oplus b_{2n})=1
\]
and
\[
 a\oplus b_{2n+1}=1,\qquad
 \kappa(a\oplus b_{2n+1})=0,
\]
 by \(j\circ\lambda_a=\widehat{\lambda_a} \circ j\),  we have
\[
 \widehat{\lambda_a}(j(b_{2n}))=j(\lambda_a(b_{2n}))=j(a\oplus b_{2n})\longrightarrow(1,1)\]
and
 \[
 \widehat{\lambda_a}(j(b_{2n+1}))=j(\lambda_a(b_{2n+1}))j(a\oplus b_{2n+1})\longrightarrow(1,0).
\]
This is contradiction. Thus the fixed translation has no continuous extension,
and neither does the binary operation.
\end{proof}

\begin{exa}\label{lem:original-boolean}
There is a second countable and metrizable topological
MV-algebra but not locally order convex, and hence, its MV-operation \(\oplus\) is not uniformly continuous for the canonical difference uniformity.
\end{exa}

\begin{proof}
This example comes from \cite[Example~1.1.11(a)]{Weber1991}. For the sake of completeness, we give the construction.

Let \(\mathbb{A}\) be the finite and cofinite
subsets of \(\N\), that is, \[\mathbb{A}=\{X\subseteq \N:X \text{~is finite or cofinite~}\}.\]
Then \(\mathbb{A}\) with the symmetric difference \(\triangle\) is a Boolean group. Indeed, \(\emptyset\) is the identity in \((\mathbb{A},\triangle)\), the inverse element of each \(X\in\mathbb{A}\) is itself and \(X\triangle Y=(X\cup Y)\setminus(X\cap Y)\) is in \(\mathbb{A}\) for each \(X,Y\in \mathbb{A}\). One can easily show that the operator \(\triangle\) is associative on \(\mathbb{A}\). Put \[
 B_k=\{2k-1,2k\}
\]
for each \(k\in \N \setminus \{0\}\) and
\begin{equation}\label{eq:boolean-Un}
 \mathbb{U}_n=\left\{\bigcup_{k\in F}B_k:
           F\subseteq\{n,n+1,\ldots\}\text{ is finite}\right\}.
\end{equation}
for each \(n\in \N\). Then each \(\mathbb{U}_n\) is a subgroup of the Boolean group
\((\mathbb{A},\triangle)\) such that \(\bigcap_{n\in \N}\mathbb{U}_n=\{\emptyset\}\) and \(\mathbb{U}_n\supseteq \mathbb{U}_m\) when \(m\geq n\), because \((\bigcup_{k\in F_1}B_k)\triangle(\bigcup_{k\in F_2}B_k)=\bigcup_{k\in F_1\triangle F_2}B_k\). Put
\begin{equation}\label{T}
\tau=\{\mathbb{O}\subseteq \mathbb{A}:\text{~each~}X\in \mathbb{O}\text{~there some~}n\in \N\text{~such that~} X\triangle \mathbb{U}_n\subseteq \mathbb{O}\}.
\end{equation}
Then the Boolean group
\((\mathbb{A},\triangle)\) with the topology \(\tau\) is a Hausdorff first countable topological group. Hence, it is second countable and metrizeble, because \(\mathbb{A}\) is countable.

Define \(\oplus\) and \({}^*\) on \(\mathbb{A}\) by
\[
 X\oplus Y=X\cup Y \text{~and~} X^*=\N\setminus X
\]
for each \(X,Y\in \mathbb{A}\), respectively. Then \((\mathbb{A},\oplus,{}^*, 0)\) is an MV-algebra, where \(0=\emptyset\). We claim that the MV-algebra \((\mathbb{A},\oplus,{}^*, 0)\) with the topology \(\tau\) defined by (\ref{T}) is a topological MV-algebra.

Since
\begin{align*}
(X\triangle \mathbb{U}_n)^*&=\{(X\triangle Y)^*:Y\in \mathbb{U}_n\}\\
&=\{\N\setminus(X\triangle Y):Y\in \mathbb{U}_n\}\\
&=\{(\N\setminus X)\triangle Y:Y\in \mathbb{U}_n\}\\
&=\{ X^*\triangle Y:Y\in \mathbb{U}_n\}\\
&=X^*\triangle \mathbb{U}_n,
\end{align*}
the operator \({}^*\) is continuous.

Next we shall prove that \(\oplus:\mathbb{A}\times \mathbb{A}\rightarrow \mathbb{A}\) is continuous. Take any \(X,Y\in \mathbb{A}\) and any neighbourhood \(V\) of \(X\oplus Y\). Then there some \(n\in \N\)
such that \((X\oplus Y)\triangle\mathbb{U}_n\subseteq V\). Since \(X, Y\) are either finite or cofinite sets, there is \(m\geq n\) such that either \(B_k \subseteq Z\) or \(B_k \cap Z=\emptyset\) for each \(k\geq m\), where \(Z=X,Y\). We claim that \[(X \triangle \mathbb{U}_m)\oplus (Y \triangle \mathbb{U}_m)\subseteq (X\oplus Y)\triangle\mathbb{U}_n\subseteq V,\] which implies that \(\oplus\) is continuous.

Indeed, take \(E,F\in \mathbb{U}_m\). Then
\begin{align*}
((X\triangle E)\oplus (Y\triangle F))\triangle(X\oplus Y)&=((X\triangle E)\cup (Y\triangle F))\triangle(X\cup Y)\\
&\subseteq ((X\triangle E)\triangle X)\cup ((Y\triangle F) \triangle Y) \\
&=E\cup F.
\end{align*}
Since either \(B_k \subseteq Z\) or \(B_k \cap Z=\emptyset\) for each \(k\geq m\) (where \(Z\in \{X,Y,E,F\}\)), we have either\[(((X\triangle E)\cup (Y\triangle F))\triangle(X\cup Y))\cap B_k=\emptyset\] or \[B_k\subseteq (((X\triangle E)\cup (Y\triangle F))\triangle(X\cup Y))\] for each \(k\geq m\). Note that \(E\cup F\) is finite. Hence, \[((X\triangle E)\cup (Y\triangle F))\triangle(X\cup Y)\in \mathbb{U}_m\subseteq \mathbb{U}_n.\]
This implies that
\[(X\triangle E)\oplus (Y\triangle F)=(X\triangle E)\cup (Y\triangle F)\in (X\cup Y)\triangle\mathbb{U}_n=(X\oplus Y)\triangle\mathbb{U}_n.\]
Thus, we have proved \[(X \triangle \mathbb{U}_m)\oplus (Y \triangle \mathbb{U}_m)\subseteq (X\oplus Y)\triangle\mathbb{U}_n\subseteq V.\]

Finally, we shall show that the topological MV-algebra \(\mathbb{A}\) is not locally order convex. If an order convex zero-neighborhood \(\mathbb{W}\subseteq \mathbb{U}_1\)
existed, it would contain some \(\mathbb{U}_n\).
For \(k\ge n\), both \(\emptyset\) and \(B_k\) would be
in \(\mathbb{W}\), forcing \(\{2k-1\}\in \mathbb{W}\), contrary to
\(\{2k-1\}\notin \mathbb{U}_1\). Thus local order convexity fails.

Hence, form Theorem \ref{thm:uniform} it follows that the MV-operation \(\oplus\) is not uniformly continuous for the canonical difference uniformity.
\end{proof}

 Let \(\mathbb{A}\) be the topological MV-algebra constructed in Example \ref{lem:original-boolean}. We retain the notation \(B_k\) and \(\mathbb{U}_n\) from the proof of Example \ref{lem:original-boolean}.
Take a point $\infty\notin\N$, and put \(\Omega=\N\cup\{\infty\}.\)
For $X\subseteq\Omega$, define the set \(P(X)\) by
\[
 P(X)=\{k\geq 1:|X\cap B_k|=1\}.
\]
Put
\begin{equation}\label{eq:B}
\mathbb{B}=\{X\subseteq\Omega:P(X)\text{ is finite}\}.
\end{equation}
Define \(\oplus\) and \({}^*\) on \(\mathbb{B}\) by
\[
 X\oplus Y=X\cup Y \text{~and~} X^*=\Omega\setminus X
\]
for each \(X,Y\in \mathbb{B}\), respectively. Then \((\mathbb{B},\oplus,{}^*, 0_{\mathbb{B}})\) is an MV-algebra, where \(0_{\mathbb{B}}=\emptyset\). Indeed, for all \(X,Y\in\mathbb{B}\),
\[
 P(\Omega\setminus X)=P(X),\qquad
 P(X\cup Y)\subseteq P(X)\cup P(Y).
\]
It follows that \(\mathbb{B}\) is closed under complementation and finite unions.
Since \(\emptyset,\Omega\in\mathbb{B}\), it is a Boolean algebra.

\begin{exa}\label{thm:boolean}
Let topological MV-algebra \(\mathbb{A}\) be from Example \ref{lem:original-boolean} and MV-algebra \(\mathbb{B}\) defined by (\ref{eq:B}). Define a map
\begin{equation}\label{eq:j}
 j:\mathbb{A}\longrightarrow\mathbb{B},\qquad
 j(X)=\begin{cases}
 X,&X\text{ is finite},\\
 X\cup\{\infty\},&X\text{ is cofinite}.
 \end{cases}
\end{equation}
Then there is a topology on MV-algebra \(\mathbb{B}\) such that \(\mathbb{B}\) is a canonically complete metrizable topological Boolean MV-algebra but not locally order convex and \(j\) is a uniform MV-embedding with dense image, for the canonical difference
uniformities. Consequently, $(\mathbb{B},j)$ realizes the canonical uniform
completion of \(\mathbb{A}\).
\end{exa}

\begin{proof}
For each $n\geq 1$, set
\begin{equation}\label{eq:Vn}
 \mathbb{V}_n=\left\{\bigcup_{k\in S}B_k:
 S\subseteq\{n,n+1,\ldots\}\right\}.
\end{equation}
Then each $\mathbb{V}_n$ is a subgroup of the Boolean group $(\mathbb{B},\triangle)$.
Indeed, if $E=\bigcup_{k\in S}B_k$ and $F=\bigcup_{k\in T}B_k$ belong
to $\mathbb{V}_n$, then
\[
 E\triangle F=\bigcup_{k\in S\triangle T}B_k\in \mathbb{V}_n.
\]
Moreover,
\[
\mathbb{V}_{n+1}\subseteq \mathbb{V}_n,\qquad
 \bigcap_{n\geq 1}\mathbb{V}_n=\{\varnothing\}.
\]
Put
\begin{equation}\label{T2}
\tau_\mathbb{B}=\{\mathbb{O}\subseteq \mathbb{B}:\text{~each~}X\in \mathbb{O}\text{~there some~}n\in \N\text{~such that~} X\triangle \mathbb{V}_n\subseteq \mathbb{O}\}.
\end{equation}
Then the Boolean group
\((\mathbb{B},\triangle)\) with the topology \(\tau_\mathbb{B}\) is a Hausdorff first countable topological group. Hence, it is  metrizeble. We claim that the MV-algebra \((\mathbb{B},\oplus,{}^*, 0_\mathbb{B})\) with the topology \(\tau_\mathbb{B}\) defined by (\ref{T2}) is a topological MV-algebra.

Complementation is continuous because, for every \(E\in \mathbb{V}_n\),
\[
 \Omega\setminus(X\triangle E)=(\Omega\setminus X)\triangle E.
\]
Thus it maps \(X\triangle \mathbb{V}_n\) onto \(X^*\triangle \mathbb{V}_n\).

To prove joint continuity of union at a fixed pair \((X,Y)\in\mathbb{B}^2\),
let \(n\geq 1\). Choose
\[
 m=\max\left\{n,\ 1+\max\bigl(P(X)\cup P(Y)\cup\{0\}\bigr)\right\}.
\]
Then neither \(|X\cap B_k|=1\) nor \(|Y\cap B_k|=1\) with \(k\geq m\).
For arbitrary \(E,F\in \mathbb{V}_m\), put
\[
 D=\bigl((X\triangle E)\cup(Y\triangle F)\bigr)
       \triangle(X\cup Y).
\]
If \(k<m\), then \(E\cap B_k=F\cap B_k=\emptyset\), so
\(D\cap B_k=\emptyset\). If \(k\geq m\), each of \(X,Y,E,F\) meets \(B_k\)
in either \(\emptyset\) or \(B_k\). Consequently, \(D\cap B_k\) is also
either \(\emptyset\) or \(B_k\). Finally, \(\infty\notin D\), because
\(E\) and \(F\) do not change membership of \(\infty\). Hence
\(D\in \mathbb{V}_m\subseteq \mathbb{V}_n\). This proves
\[
 \{X'\cup Y':X'\in X\triangle \mathbb{V}_m,\ Y'\in Y\triangle \mathbb{V}_m\}
 \subseteq (X\cup Y)\triangle \mathbb{V}_n.
\]
Thus union is jointly continuous, and \(\mathbb{B}\) is a topological Boolean
MV-algebra.

Next, we shall prove that the topological
MV-algebra \(\mathbb{B}\) is canonically complete. 
Define
\[
 \varepsilon_X=\begin{cases}1,&\infty\in X,\\0,&\infty\notin X,\end{cases}
 \qquad
 \sigma(X)=\left(\varepsilon_X,
       \{(k,X\cap B_k):k\in P(X)\}\right).
\]
Put
\[
 H_0=\emptyset,\qquad H_r=\bigcup_{k=1}^{r}B_k\quad(r\geq 1).
\]
We claim that
\begin{equation}\label{eq:pattern}
 (X,Y)\in E_{\mathbb{V}_n}\quad\Longleftrightarrow\quad
 \begin{cases}
 \sigma(X)=\sigma(Y),\\
 X\cap H_{n-1}=Y\cap H_{n-1},
 \end{cases}
\end{equation}
where \(E_{\mathbb{V}_n}=\{(X,Y)\in \mathbb{B}\times\mathbb{B}:X\setminus Y \in\mathbb{V}_n, Y\setminus X \in\mathbb{V}_n \}\) for each \(n\in \N\).

Suppose first that \((X,Y)\in E_{\mathbb{V}_n}\). Both directional differences avoid
\(\infty\) and \(H_{n-1}\), so \(X\) and \(Y\) agree there. On any block \(B_k\)
with \(k\geq n\), both directional differences must be either empty or the
whole block. If \(X\cap B_k\) is a singleton, this forces
\(Y\cap B_k=X\cap B_k\): any other choice would make at least one
directional difference a singleton on that block. The same argument with
\(X\) and \(Y\) interchanged shows that they have identical patterns.

Conversely, suppose the right-hand side of \eqref{eq:pattern} holds.
On each split block, \(X\) and \(Y\) agree. On every other block, each is
either empty or the whole block. Agreement on \(H_{n-1}\) and at \(\infty\)
then shows that both directional differences are unions of whole blocks
with indices at least \(n\). They belong to \(\mathbb{V}_n\), proving the claim.

Let \(\mathcal {F}\) be any Cauchy filter on \((\mathbb{B},\mathcal {U}_{\mathbb{B}})\).
For every $n\geq 1$, choose
\[
 Q_n\in\mathcal {F},\qquad Q_n\times Q_n\subseteq E_{\mathbb{V}_n},
\]
and put
\[
 G_n=\bigcap_{i=1}^{n}Q_i.
\]
Then
\begin{equation}\label{eq:Gn}
 G_n\in\mathcal {F},\qquad G_n\ne\emptyset,\qquad
 G_{n+1}\subseteq G_n,\qquad G_n\times G_n\subseteq E_{\mathbb{V}_n}.
\end{equation}
By \eqref{eq:pattern} with \(n=1\), all members of \(G_1\) have the same
pattern. Fix $Y_0\in G_1$, and write
\[
 \varepsilon=\varepsilon_{Y_0},\qquad F=P(Y_0).
\]
The set \(F\) is finite, and every member of every \(G_n\) has pattern
\(\sigma(Y_0)\).

For each \(k\geq 1\), the inclusion
\(G_{k+1}\times G_{k+1}\subseteq E_{\mathbb{V}_{k+1}}\) implies that all members of
\(G_{k+1}\) agree on \(B_k\). Hence the following definition is independent
of the choice of \(Y\):
\[
 T_k=Y\cap B_k\qquad(Y\in G_{k+1}).
\]
If \(k\in F\), then \(T_k\) is the singleton prescribed by \(\sigma(Y_0)\).
If \(k\notin F\), then \(T_k\) is either \(\emptyset\) or \(B_k\).
Set
\[
 I_\varepsilon=\begin{cases}
 \varnothing,&\varepsilon=0,\\
 \{\infty\},&\varepsilon=1,
 \end{cases}
 \qquad
 X=\left(\bigcup_{k\geq 1}T_k\right)\cup I_\varepsilon.
\]
Then \(P(X)=F\), so \(X\in\mathbb{B}\), and $\sigma(X)=\sigma(Y_0)$.

We show that \(\mathcal {F}\) converges to \(X\). Fix \(n\geq 1\) and take any
\(Y\in G_n\). If \(k<n\), then\(G_n\subseteq G_{k+1}\), whence
\(Y\cap B_k=T_k=X\cap B_k\). Also, \(X\) and \(Y\) have the same pattern.
Equation~\eqref{eq:pattern} therefore gives \((X,Y)\in E_{\mathbb{V}_n}\). Thus
\[
 G_n\subseteq E_{\mathbb{V}_n}[X].
\]
Since \(G_n\in\mathcal {F}\), it follows that \(E_{\mathbb{V}_n}[X]\in\mathcal {F}\) for every \(n\).
The sets \(E_{\mathbb{V}_n}[X]\) form a neighborhood base for the canonical uniform
topology, so \(\mathcal {F}\to X\). Every Cauchy filter therefore converges in
\(\mathbb{B}\), proving canonical completeness.

The map \(j\) in \eqref{eq:j} is injective because
\(j(X)\cap\N=X\). Complementation interchanges the finite and cofinite
cases, so
\[
 j(\N\setminus X)=\Omega\setminus j(X).
\]
Moreover, the union of two members of \(\mathbb{A}\) is cofinite if and only if at
least one of them is cofinite. Consequently,
\[
 j(X\cup Y)=j(X)\cup j(Y),\qquad j(\emptyset)=\emptyset.
\]
Thus \(j\) is a Boolean MV-homomorphism. 

We have the exact identity
\begin{equation}\label{eq:trace}
 \mathbb{V}_n\cap j(\mathbb{A})=j(\mathbb{U}_n)\qquad(n\geq 1).
\end{equation}
Indeed, the image of a cofinite set contains \(\infty\) and hence cannot
belong to \(\mathbb{V}_n\). The elements on the left are therefore precisely the
finite unions of whole blocks with indices at least \(n\), which are the
elements on the right. Preservation of symmetric difference now gives
\[
 (j(X)\triangle \mathbb{V}_n)\cap j(\mathbb{A})=j(X\triangle \mathbb{U}_n).
\]
Hence \(j\) is a topological embedding.

The canonical basic entourages of \(\mathbb{A}\) are
\[
 E^\mathbb{A}_{\mathbb{U}_n}=\{(X,Y)\in \mathbb{A}^2:X\setminus Y\in \mathbb{U}_n,\ Y\setminus X\in \mathbb{U}_n\}.
\]
Since \(j\) preserves directional differences, \eqref{eq:trace} yields
\[
 (j(X),j(Y))\in E_{\mathbb{V}_n}
 \quad\Longleftrightarrow\quad
 (X,Y)\in E^\mathbb{A}_{\mathbb{U}_n}.
\]
The canonical uniformity of \(\mathbb{A}\) is therefore exactly the pullback of
the restricted canonical uniformity of \(\mathbb{B}\). This proves that \(j\) is
a uniform embedding.

Finally, we shall show that \(j(\mathbb{A})\) is dense in \(\mathbb{B}\). Fix \(X\in\mathbb{B}\), and put \(r_0=\max(P(X)\cup\{1\})\). For every integer
\(r\geq r_0\), define
\[
 a_r=\begin{cases}
 X\cap H_r,&\infty\notin X,\\[1mm]
 (X\cap H_r)\cup(\N\setminus H_r),&\infty\in X.
 \end{cases}
\]
In the first case \(a_r\) is finite, and in the second it is cofinite;
thus \(a_r\in \mathbb{A}\). The sets \(j(a_r)\) and \(X\) agree on \(H_r\) and at
\(\infty\). All split blocks of \(X\) lie in \(H_r\), while neither set has
split blocks beyond \(r\). Thus \(\sigma(j(a_r))=\sigma(X)\), and
\eqref{eq:pattern} gives
\[
 (j(a_r),X)\in E_{\mathbb{V}_{r+1}}.
\]
For any fixed \(n\), this pair belongs to \(E_{\mathbb{V}_n}\) whenever
\(r\geq\max\{r_0,n-1\}\). Hence \(j(a_r)\to X\), proving that \(j(\mathbb{A})\)
is dense in \(\mathbb{B}\).

We have shown that \(\mathbb{B}\) is a complete separated uniform space and
\(j\) is a uniform embedding with dense image.

The natural order in these Boolean MV-algebras is inclusion. Given an order-convex
zero-neighborhood \(\mathbb{W}\) satisfying
\[
 \emptyset\in\mathbb{W}\subseteq \mathbb{V}_1.
\]
There would then exist \(n\) with \(\mathbb{V}_n\subseteq \mathbb{W}\). For any \(k\geq n\),
both \(\emptyset\) and \(B_k\) belong to \(\mathbb{W}\), and
\[
 \emptyset\subseteq\{2k-1\}\subseteq B_k.
\]
Order convexity would force \(\{2k-1\}\in \mathbb{W}\subseteq \mathbb{V}_1\), which is
impossible because every member of \(\mathbb{V}_1\) is a union of whole blocks.
Thus \(\mathbb{B}\) is not locally order convex.
\end{proof}



\section*{References}


\begin{thebibliography}{100}
\expandafter\ifx\csname url\endcsname\relax
  \def\url#1{\texttt{#1}}\fi
\expandafter\ifx\csname urlprefix\endcsname\relax\def\urlprefix{URL }\fi
\expandafter\ifx\csname href\endcsname\relax
  \def\href#1#2{#2} \def\path#1{#1}\fi



\bibitem{AT} A.V.~Arhangel'skii, M.G.~Tkachenko,
\newblock \emph{Topological Groups and Related Structures},
\newblock Atlantis Studies in Mathematics, vol.~1, Atlantis Press, Paris;
World Scientific Publishing Co. Pte. Ltd., Hackensack, NJ, 2008.

\bibitem{En} R. Engelking,
\newblock \emph{General Topology,}
\newblock Sigma Series in Pure Mathematics, Heldermann, Berlin, 1989.

\bibitem{Chang1958} C.C. Chang,
\emph{Algebraic analysis of many-valued logic},
Trans. Am. Math. Soc. 88 (1958) 467--490.

\bibitem{CDM2000}
R.~L.~O.~Cignoli, I.~M.~L.~D'Ottaviano, and D.~Mundici,
\emph{Algebraic Foundations of Many-Valued Reasoning},
Trends in Logic, vol.~7, Kluwer Academic Publishers, Dordrecht, 2000.

\bibitem{Gan2024}
A.~P.~Gan, W.~Luan, F.~A.~Deng, and Y.~C.~Yang,
\emph{On regularity of topological MV-algebras},
Fuzzy Sets and Systems \textbf{475} (2024), 108747.

\bibitem{Hoo1997} C. S. Hoo,
 \emph{Topological MV-algebras},
Topology Appl. 81 (1997), 103--121.

\bibitem{LiYang2026}
Q.~Li and Y.~Yang,
Locally convex topological MV-algebras and uniformization of MV-algebras,
\emph{Topology and its Applications}
\textbf{382} (2026), 109771.

\bibitem{LWY2021}
W.~Luan, H.~Weber, and Y.~Yang,
Filter topologies and topological MV-algebras,
\emph{Fuzzy Sets and Systems}
\textbf{406} (2021), 11--21.

\bibitem{Na1965}L. Nachbin,
 \emph{Topology and order},
  D. Van Nostrand Company, New York, 1965

\bibitem{Mundici1986}
D.~Mundici,
Interpretation of AF $C^*$-algebras in \L ukasiewicz sentential calculus,
\emph{Journal of Functional Analysis}
\textbf{65} (1986), 15--63.

\bibitem{RD1981}
W. Roelcke, S. Dierolf,
 Uniform structures on topological groups and their quotients,
 \emph{McGraw-Hill International Book Co., New York} 1981.

\bibitem{Weber1991}
H.~Weber,
Uniform lattices. I: A generalization of topological Riesz spaces
and topological Boolean rings,
\emph{Annali di Matematica Pura ed Applicata} (4) \textbf{160} (1991), 347--370.

\bibitem{Weber2012}
H.~Weber,
On topological MV-algebras and topological $\ell$-groups,
\emph{Topology and its Applications}
\textbf{159} (2012), 3392--3395.



\end{thebibliography}

\def\cprime{$'$}

\end{document}